\documentclass[12pt]{article}%
\usepackage{amsmath}
\usepackage{amsfonts}
\usepackage{amssymb}
\usepackage{graphicx}%
\newtheorem{theorem}{Theorem}

\newtheorem{corollary}[theorem]{Corollary}

\newtheorem{definition}[theorem]{Definition}

\newtheorem{lemma}[theorem]{Lemma}

\newtheorem{problem}[theorem]{Problem}
\newtheorem{proposition}[theorem]{Proposition}
\newtheorem{remark}[theorem]{Remark}

\newenvironment{proof}[1][Proof]{\noindent\textbf{#1.} }{\ \rule{0.5em}{0.5em}}
\begin{document}

\date{}
\title{Some remarks on a paper by L. Carlitz}
\author{Diego Dominici \thanks{e-mail: dominicd@newpaltz.edu}\\Department of Mathematics\\State University of New York at New Paltz\\75 S. Manheim Blvd. Suite 9\\New Paltz, NY 12561-2443\\USA\\Phone: (845) 257-2607\\Fax: (845) 257-3571 }
\maketitle

\begin{abstract}
We study a family of orthogonal polynomials which generalizes a sequence of
polynomials considered by L. Carlitz. We show that they are a special case of
the Sheffer polynomials and point out some interesting connections with
certain Sobolev orthogonal polynomials.

\end{abstract}


MSC-class: 33A65 (Primary) 11B83, 46E39 (Secondary)

\section{Introduction}

In \cite{MR0109229} R. Kelisky defined a set of rational integers $T_{n}$ by
means of the exponential generating function
\[
\exp\left[  \arctan\left(  z\right)  \right]  =%
{\displaystyle\sum\limits_{n=0}^{\infty}}
T_{n}\frac{z^{n}}{n!}.
\]
L. Carlitz extended Kelisky's idea by defining the polynomials $T_{n}(x),$
with%
\begin{equation}
\exp\left[  x\arctan\left(  z\right)  \right]  =%
{\displaystyle\sum\limits_{n=0}^{\infty}}
T_{n}(x)\frac{z^{n}}{n!}. \label{Tn}%
\end{equation}
In his article \cite{MR0109230}, he showed that $T_{n}(x)$ satisfies the
recurrence%
\begin{equation}
T_{n+1}=xT_{n}-n(n-1)T_{n-1},\quad n\geq0, \label{recTn}%
\end{equation}
where $T_{-1}(x)=0$ and $T_{0}(x)=1.$ He also proved \cite{MR0175839} several
arithmetic properties of $T_{n}(x)$.

In \cite{MR1742178} S. Kaijser considered the recurrence (\ref{recTn}) and
concluded that:

\begin{enumerate}
\item The function $T_{n}(x)$ is a monic polynomial of degree $n.$

\item The set $\left\{  \frac{1}{n!}T_{n}(x)\ |\ n\geq0\right\}  $ is an
orthonormal basis of the Hilbert space $H^{2}(S,P),$ where $S=\left\{
z\in\mathbb{C}\ |\ \operatorname{Im}(z)\in\lbrack-1,1]\right\}  $ and $P$ is
the Poisson measure for $0.$

\item The norm of the polynomial $\frac{1}{n!}T_{n}(x)$ is $\sqrt{2}$ if
$n\geq1$ and $1$ if $n=0.$

\item The exponential generating function of $T_{n}(x)$ is $\exp\left[
x\arctan\left(  z\right)  \right]  .$
\end{enumerate}

The same results and some extensions were presented by T. Araaya in
\cite{MR2075009}. He also showed that
\begin{equation}
\frac{1}{n!}T_{n}(x)=\frac{x}{n}P_{n-1}^{(1)}\left(  \frac{x}{2};\frac{\pi}%
{2}\right)  ,\quad n\geq1 \label{TMP}%
\end{equation}
where $P_{n}^{(\lambda)}\left(  x;\phi\right)  $ is the Meixner-Pollaczek
polynomial \cite{MR0481884}.

In this paper we extend (\ref{Tn}) by consider the following problem.

\begin{problem}
\label{problem 1} For which functions $f(z),$ will the polynomials $\Psi
_{n}(x)$, generated by%
\begin{equation}
\exp\left[  xf\left(  z\right)  \right]  =%
{\displaystyle\sum\limits_{n=0}^{\infty}}
\Psi_{n}(x)\frac{z^{n}}{n!}, \label{gen}%
\end{equation}
form an orthogonal set?
\end{problem}

In Section 2 we show that the, somewhat surprising, answer is that the
function $f(z)$ must be of the $\arctan\left(  \cdot\right)  $ type.

Problem \ref{problem 1} is a particular case of the Sheffer problem
\cite{MR0000081} of characterizing the orthogonal polynomials $S_{n}(x)$
generated by
\[
F(z)\exp\left[  xf\left(  z\right)  \right]  =%
{\displaystyle\sum\limits_{n=0}^{\infty}}
\Psi_{n}(x)\frac{z^{n}}{n!}%
\]
where $F(0)=1$ and $f(0)=0.$ Although the Sheffer polynomials have been
studied extensively \cite{MR0262569}, \cite{MR0296361}, \cite{MR1291359},
\cite{MR1284707}, \cite{MR689029}, \cite{MR1619516} \cite{MR1667138},
\cite{MR0005168}, the limiting case (\ref{gen}) with $F(z)=1$ does not appear
to have been considered before.

It is known \cite{MR1100286}, \cite{MR0481884}, that the Sheffer polynomials
reduce to the Hermite, Laguerre, Charlier, Meixner or Meixner-Pollaczek
polynomials, depending on the choice of $F(z)$ and $f\left(  z\right)  .$
Motivated by (\ref{TMP}), we wondered if some of the polynomials generated by
(\ref{gen}) will form a new class. In Section 3 we show that they will be
limiting cases of the Laguerre, Meixner-Pollaczek or Meixner polynomials.

An interesting property of these polynomials is that they do not form an
orthogonal set with respect to the standard inner-products, but they are
orthogonal with respect to some new inner-products involving differential or
difference operators. These classes of polynomials (called Sobolev orthogonal
polynomials) have been the object of much attention in the last years (see
\cite{MR1658734}, \cite{MR1013456}, \cite{MR1054761}, \cite{MR2016838},
\cite{MR1882484}, \cite{MR1643881},\ \cite{MR1104157}, \cite{MR1690024},
\cite{MR1379129}, \cite{MR1433510}, \cite{MR1608708}, \cite{MR1386733},
\cite{MR1651585}, \cite{MR1246854}, \cite{MR1891026}, \cite{MR1307190},
\cite{MR1343524}, \cite{MR1371348}, \cite{MR1393104}, \cite{MR1400044}).

\section{The polynomials $\Psi_{n}(x)$}

\begin{definition}
A function $\mu(x):R\rightarrow R$ is called a \textit{distribution function} if:

\begin{enumerate}
\item $\mu(x)$ is non-decreasing.

\item $\mu(x)$ is bounded.

\item $\mu(x)$ has finite moments, i.e.,$%
{\displaystyle\int\limits_{-\infty}^{\infty}}
x^{n}d\mu<\infty,\quad n\geq0.$
\end{enumerate}

The set
\[
\mathfrak{S}\left(  \mu\right)  =\left\{  x\in\mathbb{R}\ |\ \mu\left(
x^{+}\right)  -\mu\left(  x^{-}\right)  >0\right\}
\]
is called the spectrum of $\mu(x)$ .
\end{definition}

We remind the reader of the following result, a proof of which can be found in
\cite{MR1376370}.

\begin{theorem}
A set of orthogonal polynomials $\left\{  p_{n}(x)\ |\ n\geq0\right\}  $
satisfies a three term recurrence relation of the form%
\begin{equation}
p_{n+1}=\left(  A_{n}x+B_{n}\right)  p_{n}-C_{n}p_{n-1},\quad n\geq0
\label{rec}%
\end{equation}
with $p_{-1}(x)=0$ and $p_{0}(x)=1.$
\end{theorem}

The converse is known as Favard's Theorem \cite{Favard}.

\begin{theorem}
Let $\left\{  p_{n}(x)\ |\ n\geq0\right\}  $ be a sequence of polynomials,
satisfying (\ref{rec}) with $A_{n}\neq0$ for all $n\geq0.$ Then,

\begin{description}
\item[(a)] There exists a function $\mu(x)$ of bounded variation on
$\mathbb{(}-\infty,\infty)$ such that%
\[%
{\displaystyle\int\limits_{-\infty}^{\infty}}
p_{n}(x)p_{m}(x)d\mu=M_{n}\delta_{n,m},\quad n,m\geq0
\]
where%
\[
M_{n}=\mu_{0}\frac{A_{n}}{A_{0}}%
{\displaystyle\prod\limits_{k=1}^{n}}
C_{k}%
\]
and $\mu_{0}>0$ is a normalization constant.

\item[(b)] $\mu(x)$ can be chosen to be real valued if and only if $A_{n},$
$B_{n}$ and $C_{n}\in\mathbb{R}$ for all $n\geq0.$

\item[(c)] $\mu(x)$ is a distribution with an infinite spectrum if and only if
$A_{n},$ $B_{n}$ and $C_{n}\in\mathbb{R}$ for all $n\geq0$ and $M_{n}>0$ for
all $n\geq1.$
\end{description}
\end{theorem}

\begin{proof}
See \cite{MR0481884} for a proof and \cite{MR0184042}, \cite{MR1962997},
\cite{MR1221902}, \cite{MR1808576}, for some generalizations.
\end{proof}

Therefore, if $\left\{  \Psi_{n}(x)\ |\ n\geq0\right\}  $ is to be an
orthogonal set, the polynomials $\Psi_{n}(x)$ must satisfy a recurrence
relation of the form
\begin{equation}
\Psi_{n+1}=\left(  A_{n}x+B_{n}\right)  \Psi_{n}-C_{n}\Psi_{n-1}
\label{recpsi}%
\end{equation}
with $\Psi_{-1}(x)=0,\ \Psi_{0}(x)=1.$

Let's define the function
\begin{equation}
G(x,z)=\exp\left[  xf\left(  z\right)  \right]  . \label{G}%
\end{equation}
From (\ref{gen}) we see that in order to have $\Psi_{0}(x)=1,$ we need
$f(0)=0.$ Also, from
\begin{equation}
\frac{\partial G}{\partial z}=xf^{\prime}\left(  z\right)  G(x,z)=%
{\displaystyle\sum\limits_{n=0}^{\infty}}
\Psi_{n+1}(x)\frac{z^{n}}{n!}, \label{Gz}%
\end{equation}
we have $f^{\prime}\left(  0\right)  x=\Psi_{1}(x).$ If $\Psi_{1}(x)$ is to be
a polynomial of degree $1,$ we need $f^{\prime}\left(  0\right)  \neq0.$

Taking the $n^{th}$ derivative with respect to $z$ in (\ref{Gz}) and using
Leibniz's rule, we get%
\begin{equation}
\frac{\partial^{n+1}G}{\partial z^{n+1}}(x,z)=x%
{\displaystyle\sum\limits_{k=0}^{n}}
\binom{n}{k}f^{\left(  n+1-k\right)  }(z)\frac{\partial^{k}G}{\partial z^{k}%
}(x,z),\quad n\geq0. \label{Gn}%
\end{equation}
Setting $z=0$ in (\ref{Gn}) and using (\ref{gen}) we obtain%
\begin{align}
\Psi_{n+1}(x)  &  =x%
{\displaystyle\sum\limits_{k=0}^{n}}
\binom{n}{k}f^{\left(  n+1-k\right)  }(0)\Psi_{k}(x)\label{psin+1}\\
&  =xf^{\prime}\left(  0\right)  \Psi_{n}(x)+x%
{\displaystyle\sum\limits_{k=0}^{n-1}}
\binom{n}{k}f^{\left(  n+1-k\right)  }(0)\Psi_{k}(x).\nonumber
\end{align}
Comparing (\ref{psin+1}) with (\ref{recpsi}) we conclude that
\begin{equation}
A_{n}=f^{\prime}\left(  0\right)  \equiv a \label{An}%
\end{equation}
and%
\begin{equation}
B_{n}\Psi_{n}-C_{n}\Psi_{n-1}=x%
{\displaystyle\sum\limits_{k=0}^{n-1}}
\binom{n}{k}f^{\left(  n+1-k\right)  }(0)\Psi_{k}(x). \label{BnCn}%
\end{equation}
Using (\ref{psin+1}) in (\ref{BnCn}) we have%
\begin{align}
&  B_{n}x%
{\displaystyle\sum\limits_{k=0}^{n-1}}
\binom{n-1}{k}f^{\left(  n-k\right)  }(0)\Psi_{k}(x)-C_{n}x%
{\displaystyle\sum\limits_{k=0}^{n-2}}
\binom{n-2}{k}f^{\left(  n-k-1\right)  }(0)\Psi_{k}(x)\label{BnCn1}\\
&  =x%
{\displaystyle\sum\limits_{k=0}^{n-1}}
\binom{n}{k}f^{\left(  n+1-k\right)  }(0)\Psi_{k}(x),\nonumber
\end{align}
which implies%
\[
B_{n}f^{\prime}\left(  0\right)  \Psi_{n-1}(x)=nf^{\prime\prime}\left(
0\right)  \Psi_{n-1}(x)
\]
or%
\begin{equation}
B_{n}=n\frac{f^{\prime\prime}\left(  0\right)  }{f^{\prime}\left(  0\right)
}\equiv nb. \label{Bn}%
\end{equation}

Rearranging terms in (\ref{BnCn1}) and using (\ref{Bn}), we get%
\[%
{\displaystyle\sum\limits_{k=0}^{n-2}}
\left[  nb\binom{n-1}{k}f^{\left(  n-k\right)  }(0)-C_{n}\binom{n-2}%
{k}f^{\left(  n-k-1\right)  }(0)-\binom{n}{k}f^{\left(  n+1-k\right)
}(0)\right]  \Psi_{k}(x)=0
\]
which gives%
\begin{equation}
\frac{-C_{n}}{n\left(  n-1\right)  }=\frac{f^{\left(  k+1\right)
}(0)-kbf^{\left(  k\right)  }(0)}{k\left(  k-1\right)  f^{\left(  k-1\right)
}(0)}. \label{sep}%
\end{equation}
For equation (\ref{sep}) to be valid for $2\leq k\leq n,\ n\geq2,$ it is
necessary that both sides be equal to the same constant $-c,$ with $c>0.$
Thus,%
\[
\frac{-C_{n}}{n\left(  n-1\right)  }=-c
\]
or%
\begin{equation}
C_{n}=cn\left(  n-1\right)  ,\quad c>0 \label{Cn}%
\end{equation}
and%
\begin{equation}
f^{\left(  k+1\right)  }(0)-kbf^{\left(  k\right)  }(0)+ck\left(  k-1\right)
f^{\left(  k-1\right)  }(0)=0,\quad k\geq1. \label{reqf}%
\end{equation}

The solution of the recurrence (\ref{reqf}), with $f^{\left(  1\right)
}(0)=a,$ is given by%
\begin{equation}
f^{\left(  k\right)  }(0)=\frac{ac}{R^{-}\left(  bR^{+}-2c\right)  }%
\Gamma(k)\left[  \left(  R^{+}\right)  ^{k}-\left(  R^{-}\right)  ^{k}\right]
,\quad k\geq1 \label{df}%
\end{equation}
where $\Gamma\left(  \cdot\right)  $ is the Gamma function and%
\begin{equation}
R^{\pm}=\frac{1}{2}\left[  b\pm\sqrt{b^{2}-4c}\right]  . \label{R}%
\end{equation}
From (\ref{df}) and $f(0)=0,$ we get%
\begin{equation}
f(z)=\frac{a}{\sqrt{b^{2}-4c}}\ln\left[  \frac{1-zR^{-}}{1-zR^{+}}\right]
\label{f}%
\end{equation}
or%
\[
f(z)=\frac{2a}{\sqrt{b^{2}-4c}}\left[  \operatorname{arctanh}\left(
\frac{b-2cz}{\sqrt{b^{2}-4c}}\right)  -\operatorname{arctanh}\left(  \frac
{b}{\sqrt{b^{2}-4c}}\right)  \right]  .
\]

We summarize the results of this section in the following theorem.

\begin{theorem}
If the family of polynomials $\left\{  \Psi_{n}(x)\ |\ n\geq0\right\}  $
defined by (\ref{gen}) satisfies (\ref{recpsi}), then%
\begin{equation}
f(z)=\frac{2a}{\sqrt{b^{2}-4c}}\left[  \operatorname{arctanh}\left(
\frac{b-2cz}{\sqrt{b^{2}-4c}}\right)  -\operatorname{arctanh}\left(  \frac
{b}{\sqrt{b^{2}-4c}}\right)  \right]  , \label{tanh}%
\end{equation}
and%
\[
A_{n}=a,\quad B_{n}=nb,\quad C_{n}=cn\left(  n-1\right)
\]
with $a\neq0$ and $c>0.$ The recurrence relation (\ref{recpsi}) takes the
form
\begin{equation}
\Psi_{n+1}=\left(  ax+nb\right)  \Psi_{n}-cn\left(  n-1\right)  \Psi
_{n-1},\quad n\geq0. \label{recpsi1}%
\end{equation}

\end{theorem}

\begin{remark}
In the remainder of the paper we shall stress the dependence of $\Psi_{n}(x)$
on $a,b$ and $c$ by writing $\Psi_{n}(x)=\Psi_{n}(x;a,b,c).$ Note that from
(\ref{tanh}) it follows that
\begin{equation}
T_{n}(x)=\Psi_{n}(x;1,0,1).\label{TPsi}%
\end{equation}

\end{remark}

\begin{remark}
Since $C_{1}=0,$ it follows from Favard's theorem that the full set $\left\{
\Psi_{n}(x)\ |\ n\geq0\right\}  $ is not orthogonal with respect to an
inner-product generated by a distribution function. However, we shall see in
the next section that $\left\{  \Psi_{n}(x)\ |\ n\geq1\right\}  $ is$.$
\end{remark}

\section{Representation of $\Psi_{n}$}

In order to find a representation of $\Psi_{n}(x;a,b,c),$ in terms of the
classic hypergeometric polynomials, we shall consider three different cases,
depending on the sign of the discriminant in (\ref{R}).

\subsection{Case 1: $b^{2}-4c=0,\ b\neq0$}

From (\ref{R}), we have $R^{+}=R^{-}=\frac{1}{2}b.$ Taking the limit as
$c\rightarrow\frac{b^{4}}{4}$ in (\ref{f}) gives%
\begin{equation}
f(z)=-\frac{2a}{b}\frac{\frac{b}{2}z}{\frac{b}{2}z-1}. \label{f1}%
\end{equation}
The generating function for the Laguerre polynomials is given by
\cite{MR0481145}%
\begin{equation}
\left(  1-z\right)  ^{-a-1}\exp\left[  x\frac{z}{z-1}\right]  =%
{\displaystyle\sum\limits_{n=0}^{\infty}}
L_{n}^{(\alpha)}(x)z^{n},\quad\alpha>-1. \label{La}%
\end{equation}
Hence, from (\ref{f1}) and (\ref{La}), we obtain%
\begin{equation}
\Psi_{n}\left(  x;a,b,\frac{b^{4}}{4}\right)  =n!\left(  \frac{b}{2}\right)
^{n}\underset{\alpha\rightarrow-1}{\lim}L_{n}^{(\alpha)}\left(  -\frac{2a}%
{b}x\right)  . \label{lim1}%
\end{equation}

To find the limit in (\ref{lim1}) we prove the following lemma.

\begin{lemma}
\label{lemma1} Let $_{r}F_{s}$ be the hypergeometric function \cite{MR0350075}%
. Then,
\begin{align}
&  \underset{\omega\rightarrow0}{\lim}\left(  \omega\right)  _{n}\ _{r}%
F_{s}\left(  \left.
\begin{array}
[c]{c}%
a_{1}\left(  \omega\right)  ,\ldots,a_{r}\left(  \omega\right) \\
\omega,b_{2}\left(  \omega\right)  ,\ldots,b_{s}\left(  \omega\right)
\end{array}
\right\vert x\right) \label{lemma}\\
&  =x\Gamma\left(  n\right)  \frac{%
{\displaystyle\prod\limits_{j=1}^{r}}
a_{j}\left(  0\right)  }{%
{\displaystyle\prod\limits_{j=2}^{s}}
b_{j}\left(  0\right)  }\ _{r}F_{s}\left(  \left.
\begin{array}
[c]{c}%
a_{1}\left(  0\right)  +1,\ldots,a_{r}\left(  0\right)  +1\\
2,b_{2}\left(  0\right)  +1,\ldots,b_{s}\left(  0\right)  +1
\end{array}
\right\vert x\right) \nonumber
\end{align}

\end{lemma}

\begin{proof}%
\[
\underset{\omega\rightarrow0}{\lim}\left(  \omega\right)  _{n}\ _{r}%
F_{s}\left(  \left.
\begin{array}
[c]{c}%
a_{1}\left(  \omega\right)  ,\ldots,a_{r}\left(  \omega\right) \\
\omega,b_{2}\left(  \omega\right)  ,\ldots,b_{s}\left(  \omega\right)
\end{array}
\right\vert x\right)
\]%
\[
=\underset{\omega\rightarrow0}{\lim}\left(  \omega\right)  _{n}%
{\displaystyle\sum\limits_{k=0}^{\infty}}
\frac{\left(  a_{1}\left(  \omega\right)  \right)  _{k}\times\cdots
\times\left(  a_{r}\left(  \omega\right)  \right)  _{k}}{\left(
\omega\right)  _{k}\left(  b_{2}\left(  \omega\right)  \right)  _{k}%
\times\cdots\times\left(  b_{s}\left(  \omega\right)  \right)  _{k}\left(
1\right)  _{k}}x^{k}%
\]%
\[
=\underset{\omega\rightarrow0}{\lim}\Gamma\left(  \omega+n\right)
{\displaystyle\sum\limits_{k=0}^{\infty}}
\frac{\left(  a_{1}\left(  \omega\right)  \right)  _{k}\times\cdots
\times\left(  a_{r}\left(  \omega\right)  \right)  _{k}}{\Gamma\left(
\omega+k\right)  \left(  b_{2}\left(  \omega\right)  \right)  _{k}\times
\cdots\times\left(  b_{s}\left(  \omega\right)  \right)  _{k}\left(  1\right)
_{k}}x^{k}%
\]%
\[
=\Gamma\left(  n\right)
{\displaystyle\sum\limits_{k=1}^{\infty}}
\frac{\left(  a_{1}\left(  0\right)  \right)  _{k}\times\cdots\times\left(
a_{r}\left(  0\right)  \right)  _{k}}{\Gamma\left(  k\right)  \left(
b_{2}\left(  0\right)  \right)  _{k}\times\cdots\times\left(  b_{s}\left(
0\right)  \right)  _{k}\left(  1\right)  _{k}}x^{k}%
\]%
\[
=\Gamma\left(  n\right)
{\displaystyle\sum\limits_{k=0}^{\infty}}
\frac{\left(  a_{1}\left(  0\right)  \right)  _{k+1}\times\cdots\times\left(
a_{r}\left(  0\right)  \right)  _{k+1}}{\Gamma\left(  k+1\right)  \left(
b_{2}\left(  0\right)  \right)  _{k+1}\times\cdots\times\left(  b_{s}\left(
0\right)  \right)  _{k+1}\left(  1\right)  _{k+1}}x^{k+1}%
\]%
\[
=x\Gamma\left(  n\right)  \frac{%
{\displaystyle\prod\limits_{j=1}^{r}}
a_{j}\left(  0\right)  }{%
{\displaystyle\prod\limits_{j=2}^{s}}
b_{j}\left(  0\right)  }%
{\displaystyle\sum\limits_{k=0}^{\infty}}
\frac{\left(  a_{1}\left(  0\right)  +1\right)  _{k}\times\cdots\times\left(
a_{r}\left(  0\right)  +1\right)  _{k}}{\left(  1\right)  _{k}\left(
b_{2}\left(  0\right)  +1\right)  _{k}\times\cdots\times\left(  b_{s}\left(
0\right)  +1\right)  _{k}\left(  2\right)  _{k}}x^{k}%
\]%
\[
=x\Gamma\left(  n\right)  \frac{%
{\displaystyle\prod\limits_{j=1}^{r}}
a_{j}\left(  0\right)  }{%
{\displaystyle\prod\limits_{j=2}^{s}}
b_{j}\left(  0\right)  }\ _{r}F_{s}\left(  \left.
\begin{array}
[c]{c}%
a_{1}\left(  0\right)  +1,\ldots,a_{r}\left(  0\right)  +1\\
2,b_{2}\left(  0\right)  +1,\ldots,b_{s}\left(  0\right)  +1
\end{array}
\right\vert x\right)
\]

\end{proof}

\begin{corollary}
We can extend the class of Laguerre polynomials by defining
\begin{equation}
L_{0}^{(-1)}=1,\quad L_{n}^{(-1)}\left(  x\right)  =-\frac{x}{n}L_{n-1}%
^{(1)}\left(  x\right)  ,\quad n\geq1. \label{L-1}%
\end{equation}

\end{corollary}

\begin{proof}
Using the definition \cite{koekoek94askeyscheme}%
\begin{equation}
L_{n}^{(\alpha)}\left(  x\right)  =\frac{\left(  \alpha+1\right)  _{n}}%
{n!}\ _{1}F_{1}\left(  \left.
\begin{array}
[c]{c}%
-n\\
\alpha+1
\end{array}
\right\vert x\right)  , \label{Ldef}%
\end{equation}
and (\ref{lemma}), we get%
\[
\underset{\alpha\rightarrow-1}{\lim}L_{n}^{(\alpha)}\left(  x\right)
=\frac{1}{n!}x\Gamma\left(  n\right)  \left(  -n\right)  \ _{1}F_{1}\left(
\left.
\begin{array}
[c]{c}%
-n+1\\
2
\end{array}
\right\vert x\right)  .
\]
From (\ref{Ldef}) we have%
\[
_{1}F_{1}\left(  \left.
\begin{array}
[c]{c}%
-n+1\\
2
\end{array}
\right\vert x\right)  =\frac{\left(  n-1\right)  !}{\left(  2\right)  _{n-1}%
}L_{n-1}^{(1)}\left(  x\right)  =\frac{1}{n}L_{n-1}^{(1)}\left(  x\right)
,\quad n\geq1.
\]
Therefore,%
\[
\underset{\alpha\rightarrow-1}{\lim}L_{n}^{(\alpha)}\left(  x\right)
=-\frac{x}{n}L_{n-1}^{(1)}\left(  x\right)  ,\quad n\geq1.
\]

\end{proof}

\begin{corollary}
We have the representation
\[
\Psi_{n}\left(  x;a,b,\frac{b^{4}}{4}\right)  =a\left(  \frac{b}{2}\right)
^{n-1}\left(  n-1\right)  !L_{n-1}^{(1)}\left(  -\frac{2a}{b}x\right)  ,\quad
n\geq1.
\]

\end{corollary}

\begin{remark}
If we define the inner product%
\begin{equation}
\left\langle f,g\right\rangle _{\alpha}=%
{\displaystyle\int\limits_{-\infty}^{\infty}}
f(x)g(x)e^{-x}x^{\alpha}\chi_{\lbrack0,\infty)}(x)dx, \label{Ip1}%
\end{equation}
\medskip we have \cite{koekoek94askeyscheme}%
\begin{gather*}
\left\langle L_{n}^{(-1)},L_{m}^{\left(  -1\right)  }\right\rangle
_{-1}=\left\langle -\frac{x}{n}L_{n-1}^{(1)},-\frac{x}{m}L_{m-1}%
^{(1)}\right\rangle _{-1}\\
=\frac{1}{nm}%
{\displaystyle\int\limits_{0}^{\infty}}
L_{n-1}^{(1)}(x)L_{m-1}^{(1)}(x)x^{2}e^{-x}x^{-1}dx=\frac{1}{nm}\left\langle
L_{n-1}^{(1)},L_{m-1}^{(1)}\right\rangle _{1}\\
=\frac{1}{n^{2}}\frac{\Gamma\left(  n+1\right)  }{\left(  n-1\right)  !}%
\delta_{n-1,m-1}=\frac{1}{n}\delta_{n,m},\quad n,m\geq1.
\end{gather*}
Therefore, $\left\{  L_{n}^{(-1)}\ |\ n\geq1\right\}  $is an orthogonal set
with respect to the inner product $\left\langle \cdot,\cdot\right\rangle
_{-1}$ defined by (\ref{Ip1}). The general case was studied in
\cite{MR2077205} and \cite{MR1333756}, where it was shown that $\left\{
L_{n}^{(-k)}\ |\ n\geq k\right\}  $is an orthogonal set with respect to the
inner product $\left\langle \cdot,\cdot\right\rangle _{-k}$ for all $k\geq1.$
\end{remark}

\begin{remark}
We will now show that $\left\{  L_{n}^{(-1)}\ |\ n\geq0\right\}  $is an
orthonormal set with respect to the inner product%
\[
\left\langle f,g\right\rangle =f(0)g(0)+%
{\displaystyle\int\limits_{0}^{\infty}}
f^{\prime}(x)g^{\prime}(x)e^{-x}dx.
\]

From (\ref{L-1}) we have%
\begin{equation}
L_{n}^{(-1)}(0)=\delta_{n,0},\quad n\geq0, \label{L-10}%
\end{equation}
which implies%
\[
\left\langle L_{0}^{(-1)},L_{n}^{(-1)}\right\rangle =L_{n}^{(-1)}(0)+%
{\displaystyle\int\limits_{0}^{\infty}}
0e^{-x}dx=\delta_{n,0},\quad n\geq0.
\]
Using the formula \cite{MR0350075}
\[
\frac{d}{dx}L_{n}^{(\alpha)}=-L_{n-1}^{(\alpha+1)},\quad n\geq0
\]
and (\ref{L-10}), we get%
\begin{align}
\left\langle L_{n}^{(-1)},L_{m}^{(-1)}\right\rangle  &  =L_{n}^{(-1)}%
(0)L_{m}^{(-1)}(0)+%
{\displaystyle\int\limits_{0}^{\infty}}
L_{n-1}^{(0)}(x)L_{m-1}^{(0)}(x)e^{-x}dx\label{L-1n}\\
&  =\delta_{n-1,m-1}=\delta_{n,m},\quad n,m\geq1.\nonumber
\end{align}
Thus, (\ref{L-10}) and (\ref{L-1n}) give%
\[
\left\langle L_{n}^{(-1)},L_{m}^{(-1)}\right\rangle =\delta_{n,m},\quad
n,m\geq0.
\]

The general case for $L_{n}^{(-k)}$ was first considered in \cite{MR1343537}
and subsequently in \cite{MR1433121}, \cite{MR1488822}, \cite{MR1399895},
\cite{MR1991819}, \cite{MR1213122}, \cite{MR1954934} and \cite{MR1662711}.
\end{remark}

\subsection{Case 2: $b^{2}-4c<0$}

From (\ref{R}) we see that in this case $R^{\pm}$ are complex conjugates,
\begin{equation}
R^{\pm}=\frac{1}{2}\left[  b\pm\mathrm{i}\sqrt{4c-b^{2}}\right]  \label{R1}%
\end{equation}
with absolute value $\left\vert R^{\pm}\right\vert =\sqrt{c}.$ Hence, we have%
\begin{equation}
\frac{R^{\pm}}{\sqrt{c}}=e^{\mathrm{i}\phi},\quad0<\phi<\pi. \label{phi}%
\end{equation}
Using (\ref{phi}) in (\ref{f}) we can write $f(z)$ as%
\begin{equation}
f(z)=\frac{a\mathrm{i}}{\sqrt{4c-b^{2}}}\ln\left[  \frac{1-z\sqrt
{c}e^{\mathrm{i}\phi}}{1-z\sqrt{c}e^{-\mathrm{i}\phi}}\right]  . \label{f2}%
\end{equation}

Comparing (\ref{f2}) and the generating function for the Meixner-Pollaczek
polynomials, given by
\begin{equation}
\left(  1-ze^{\mathrm{i}\phi}\right)  ^{-\lambda+\mathrm{i}x}\left(
1-ze^{-\mathrm{i}\phi}\right)  ^{-\lambda-\mathrm{i}x}=%
{\displaystyle\sum\limits_{n=0}^{\infty}}
P_{n}^{(\lambda)}(x;\phi)z^{n},\quad\lambda>0,\quad0<\phi<\pi,\label{MP}%
\end{equation}
we conclude that%
\begin{equation}
\Psi_{n}\left(  x;a,b,c\right)  =n!c^{\frac{n}{2}}\underset{\lambda
\rightarrow0}{\lim}P_{n}^{(\lambda)}\left(  \frac{ax}{\sqrt{4c-b^{2}}}%
;\phi\right)  ,\quad b^{2}-4c<0\label{psi2}%
\end{equation}
with $\phi$ defined by (\ref{phi}). We will then find the limit in
(\ref{psi2}) using Lemma \ref{lemma1}.

\begin{proposition}
The family of Meixner-Pollaczek polynomials can be extended if we define%
\begin{equation}
P_{0}^{(0)}\left(  x;\phi\right)  =1,\quad P_{n}^{(0)}\left(  x;\phi\right)
=\frac{2x}{n}\sin\left(  \phi\right)  P_{n-1}^{(1)}\left(  x;\phi\right)
,\quad n\geq1. \label{MP0}%
\end{equation}

\end{proposition}

\begin{proof}
Using the definition \cite{koekoek94askeyscheme}%
\begin{equation}
P_{n}^{(\lambda)}(x;\phi)=\frac{\left(  2\lambda\right)  _{n}}{n!}%
\ e^{\mathrm{i}n\phi}\ _{2}F_{1}\left(  \left.
\begin{array}
[c]{c}%
-n,\ \lambda+\mathrm{i}x\\
2\lambda
\end{array}
\right\vert 1-e^{-2\mathrm{i}\phi}\right)  \label{Mp}%
\end{equation}
and (\ref{lemma}), we get%
\begin{gather*}
\underset{\lambda\rightarrow0}{\lim}P_{n}^{(\lambda)}(x;\phi)=\frac{1}%
{n!}\ e^{\mathrm{i}n\phi}\left(  1-e^{-2\mathrm{i}\phi}\right)  \Gamma\left(
n\right)  \left(  -n\right)  \left(  \mathrm{i}x\right)  \ \\
\times\ _{2}F_{1}\left(  \left.
\begin{array}
[c]{c}%
-n+1,\ \mathrm{i}x+1\\
2
\end{array}
\right\vert 1-e^{-2\mathrm{i}\phi}\right) \\
=\frac{1}{n!}\ e^{\mathrm{i}n\phi}\left(  1-e^{-2\mathrm{i}\phi}\right)
\Gamma\left(  n\right)  \left(  -n\right)  \left(  \mathrm{i}x\right)
P_{n-1}^{(1)}(x;\phi)\frac{(n-1)!}{\left(  2\right)  _{n-1}}e^{-\mathrm{i}%
\left(  n-1\right)  \phi}\\
=\frac{1}{n}\ \left(  \frac{e^{\mathrm{i}\phi}-e^{-\mathrm{i}\phi}}%
{\mathrm{i}}\right)  xP_{n-1}^{(1)}(x;\phi),\quad n\geq1
\end{gather*}
and the result follows.
\end{proof}

\begin{corollary}
If $b^{2}-4c<0,$ we have the representation
\begin{equation}
\Psi_{n}\left(  x;a,b,c\right)  =\left(  n-1\right)  !c^{\frac{n-1}{2}%
}axP_{n-1}^{(1)}\left(  \frac{a}{\sqrt{4c-b^{2}}}x;\phi\right)  ,\quad n\geq1
\label{psi3}%
\end{equation}
with $\phi$ defined by (\ref{phi}).
\end{corollary}

\begin{proof}
Using (\ref{MP0}) in (\ref{psi2}) we have%
\[
\Psi_{n}\left(  x;a,b,c\right)  =n!c^{\frac{n}{2}}\frac{2x}{n}\frac{a}%
{\sqrt{4c-b^{2}}}\sin\left(  \phi\right)  P_{n-1}^{(1)}\left(  \frac{a}%
{\sqrt{4c-b^{2}}}x;\phi\right)  ,\quad n\geq1.
\]
From (\ref{R1}) and (\ref{phi}) we have%
\[
\frac{1}{2\sqrt{c}}\left[  b\pm\mathrm{i}\sqrt{4c-b^{2}}\right]  =\cos\left(
\phi\right)  \pm\mathrm{i}\sin\left(  \phi\right)  ,
\]
which gives%
\[
\sin\left(  \phi\right)  =\sqrt{1-\frac{b^{2}}{4c}}.
\]

\end{proof}

\begin{remark}
Replacing $a=1,b=0$ and $c=1$ in (\ref{psi3}) and using (\ref{TPsi}), we
obtain%
\[
T_{n}(x)=\Psi_{n}(x;1,0,1)=\left(  n-1\right)  !xP_{n-1}^{(1)}\left(  \frac
{x}{2};\frac{\pi}{2}\right)  ,\quad n\geq1
\]
which agrees with (\ref{TMP}).
\end{remark}

\begin{remark}
If we define the inner product%
\begin{equation}
\left\langle f,g\right\rangle _{\lambda}=\frac{1}{2\pi}%
{\displaystyle\int\limits_{-\infty}^{\infty}}
f(x)g(x)e^{\left(  2\phi-\pi\right)  x}\left\vert \Gamma\left(  \lambda
+\mathrm{i}x\right)  \right\vert ^{2}dx, \label{IPm}%
\end{equation}
\medskip we have \cite{koekoek94askeyscheme}%
\begin{gather*}
\left\langle P_{n}^{(0)},P_{m}^{(0)}\right\rangle _{0}=\left\langle \frac
{2x}{n}\sin\left(  \phi\right)  P_{n-1}^{(1)}\left(  x;\phi\right)  ,\frac
{2x}{m}\sin\left(  \phi\right)  P_{m-1}^{(1)}\left(  x;\phi\right)
\right\rangle _{0}\\
=\frac{4\sin^{2}\left(  \phi\right)  }{nm}\frac{1}{2\pi}%
{\displaystyle\int\limits_{-\infty}^{\infty}}
P_{n-1}^{(1)}\left(  x;\phi\right)  P_{m-1}^{(1)}\left(  x;\phi\right)
e^{\left(  2\phi-\pi\right)  x}\left\vert \Gamma\left(  \mathrm{i}x\right)
\right\vert ^{2}x^{2}dx\\
=\frac{4\sin^{2}\left(  \phi\right)  }{nm}\frac{1}{2\pi}%
{\displaystyle\int\limits_{-\infty}^{\infty}}
P_{n-1}^{(1)}\left(  x;\phi\right)  P_{m-1}^{(1)}\left(  x;\phi\right)
e^{\left(  2\phi-\pi\right)  x}\left\vert \Gamma\left(  1+\mathrm{i}x\right)
\right\vert ^{2}dx\\
=\frac{4\sin^{2}\left(  \phi\right)  }{nm}\left\langle P_{n-1}^{(1)}%
,P_{m-1}^{(1)}\right\rangle _{1}\\
=\frac{4\sin^{2}\left(  \phi\right)  }{n^{2}}\frac{\Gamma\left(  n+1\right)
}{4\sin^{2}\left(  \phi\right)  \left(  n-1\right)  !}\delta_{n-1,m-1}%
=\frac{1}{n}\delta_{n,m},\quad n,m\geq1.
\end{gather*}
Therefore, $\left\{  P_{n}^{(0)}\ |\ n\geq1\right\}  $is an orthogonal set
with respect to the inner product $\left\langle \cdot,\cdot\right\rangle _{0}$
defined by (\ref{IPm}).
\end{remark}

\begin{remark}
We will now show that $\left\{  P_{n}^{(0)}\ |\ n\geq0\right\}  $is an
orthonormal set with respect to the inner product%
\[
\left\langle f,g\right\rangle =f(0)g(0)+\frac{1}{4\pi\sin\left(  \phi\right)
}%
{\displaystyle\int\limits_{-\infty}^{\infty}}
{\Large \delta}f(x){\Large \delta}g(x)e^{\left(  2\phi-\pi\right)
x}\left\vert \Gamma\left(  \frac{1}{2}+\mathrm{i}x\right)  \right\vert
^{2}dx,
\]
where%
\[
{\Large \delta}f(x)=f\left(  x+\frac{1}{2}\right)  -f\left(  x-\frac{1}%
{2}\right)  .
\]

From (\ref{MP0}) we have%
\begin{equation}
P_{n}^{(0)}(0;\phi)=\delta_{n,0},\quad n\geq0, \label{Po0}%
\end{equation}
which implies%
\begin{equation}
\left\langle P_{0}^{(0)},P_{n}^{(0)}\right\rangle =P_{n}^{(0)}(0;\phi)+%
{\displaystyle\int\limits_{-\infty}^{\infty}}
0e^{\left(  2\phi-\pi\right)  x}\left\vert \Gamma\left(  \frac{1}%
{2}+\mathrm{i}x\right)  \right\vert ^{2}dx=\delta_{n,0},\quad n\geq0.
\label{Mp0}%
\end{equation}
Using the formula \cite{koekoek94askeyscheme}
\[
{\Large \delta}P_{n}^{(\lambda)}(x;\phi)=2\sin\left(  \phi\right)
P_{n-1}^{\left(  \lambda+\frac{1}{2}\right)  }(x;\phi),\quad n\geq0
\]
and (\ref{Po0}), we get%
\begin{gather}
\left\langle P_{n}^{(0)},P_{m}^{(0)}\right\rangle =P_{n}^{(0)}(0;\phi
)P_{m}^{(0)}(0;\phi)\nonumber\\
+2\sin\left(  \phi\right)  \frac{1}{2\pi}%
{\displaystyle\int\limits_{-\infty}^{\infty}}
P_{n-1}^{\left(  \frac{1}{2}\right)  }(x;\phi)P_{m-1}^{\left(  \frac{1}%
{2}\right)  }(x;\phi)e^{\left(  2\phi-\pi\right)  x}\left\vert \Gamma\left(
\frac{1}{2}+\mathrm{i}x\right)  \right\vert ^{2}dx\label{MPn}\\
=2\sin\left(  \phi\right)  \frac{\Gamma\left(  n\right)  }{2\sin\left(
\phi\right)  \left(  n-1\right)  !}\delta_{n-1,m-1}=\delta_{n,m},\quad
n,m\geq1.\nonumber
\end{gather}
Thus, (\ref{Mp0}) and (\ref{MPn}) give%
\[
\left\langle P_{n}^{(0)},P_{m}^{(0)}\right\rangle =\delta_{n,m},\quad
n,m\geq0.
\]

\end{remark}

\subsection{Case 3: $b^{2}-4c>0$}

From (\ref{R}) we have $R^{-}<R^{+},$ if $b^{2}-4c>0.$ Since $\sqrt{b^{2}%
-4c}<\left\vert b\right\vert ,$
\[
0<R^{-}<R^{+}\text{ for \ }b>0\text{ \ and \ }R^{-}<R^{+}<0\text{ for \ }b<0.
\]
Thus,%
\begin{equation}
0<\frac{R^{-}}{R^{+}}<1\text{ for \ }b>0\text{ \ and \ }0<\frac{R^{+}}{R^{-}%
}<1\text{ for \ }b<0. \label{RP/Rm}%
\end{equation}

Comparing (\ref{f}) and the generating function for the Meixner polynomials,
given by \cite{MR1149380}
\[
\left(  1-\frac{t}{\gamma}\right)  ^{x}\left(  1-t\right)  ^{-\beta-x}=%
{\displaystyle\sum\limits_{n=0}^{\infty}}
M_{n}^{\left(  \beta\right)  }(x;\gamma)\frac{t^{n}}{n!},\quad\beta
>0,\quad0<\gamma<1,
\]
we have
\begin{equation}
\Psi_{n}\left(  x;a,b,c\right)  =\left\{
\begin{array}
[c]{c}%
\left(  R^{-}\right)  ^{n}\underset{\beta\rightarrow0}{\lim}M_{n}^{\left(
\beta\right)  }\left(  -\frac{ax}{\sqrt{b^{2}-4c}};\frac{R^{-}}{R^{+}}\right)
,\quad b>0\\
\left(  R^{+}\right)  ^{n}\underset{\beta\rightarrow0}{\lim}M_{n}^{\left(
\beta\right)  }\left(  \frac{ax}{\sqrt{b^{2}-4c}};\frac{R^{+}}{R^{-}}\right)
,\quad b<0
\end{array}
\right.  \label{psi4}%
\end{equation}
where we have used (\ref{RP/Rm}). To find the limit in (\ref{psi4}) we use
Lemma 1.

\begin{proposition}
The family of Meixner polynomials can be extended if we define%
\begin{equation}
M_{0}^{\left(  0\right)  }(x;\gamma)=1,\quad M_{n}^{\left(  0\right)
}(x;\gamma)=\left(  1-\frac{1}{\gamma}\right)  \ xM_{n-1}^{\left(  2\right)
}(x-1;\gamma),\quad n\geq1. \label{M0}%
\end{equation}

\end{proposition}

\begin{proof}
Using the definition \cite{MR1149380}%
\[
M_{n}^{\left(  \beta\right)  }(x;\gamma)=\left(  \beta\right)  _{n}\ _{2}%
F_{1}\left(  \left.
\begin{array}
[c]{c}%
-n,\ -x\\
\beta
\end{array}
\right\vert 1-\frac{1}{\gamma}\right)
\]
and (\ref{lemma}), we get%
\begin{gather*}
\underset{\beta\rightarrow0}{\lim}M_{n}^{\left(  \beta\right)  }%
(x;\gamma)=\left(  1-\frac{1}{\gamma}\right)  \ \Gamma\left(  n\right)
(-n)(-x)\ _{2}F_{1}\left(  \left.
\begin{array}
[c]{c}%
-n+1,\ -x+1\\
2
\end{array}
\right\vert 1-\frac{1}{\gamma}\right) \\
=\left(  1-\frac{1}{\gamma}\right)  \ xM_{n-1}^{\left(  2\right)  }%
(x-1;\gamma),\quad n\geq1
\end{gather*}
and the result follows.
\end{proof}

\begin{corollary}
If $b^{2}-4c>0,$ we have the representation
\begin{equation}
\Psi_{n}\left(  x;a,b,c\right)  =\left\{
\begin{array}
[c]{c}%
ax\left(  R^{-}\right)  ^{n-1}M_{n-1}^{\left(  2\right)  }\left(  -\frac
{ax}{\sqrt{b^{2}-4c}}-1;\frac{R^{-}}{R^{+}}\right)  ,\quad b>0\\
ax\left(  R^{+}\right)  ^{n-1}M_{n-1}^{\left(  2\right)  }\left(  \frac
{ax}{\sqrt{b^{2}-4c}}-1;\frac{R^{+}}{R^{-}}\right)  ,\quad b<0
\end{array}
\right.  ,\quad n\geq1. \label{psi5}%
\end{equation}

\end{corollary}

\begin{proof}
Using (\ref{M0}) in (\ref{psi5}), we have%
\begin{align*}
\Psi_{n}\left(  x;a,b,c\right)   &  =\left(  R^{-}\right)  ^{n}\left(
1-\frac{R^{-}}{R^{+}}\right)  \ \left(  -\frac{ax}{\sqrt{b^{2}-4c}}\right)
M_{n-1}^{\left(  2\right)  }\left(  -\frac{ax}{\sqrt{b^{2}-4c}}-1;\frac{R^{-}%
}{R^{+}}\right) \\
&  =\left(  R^{-}\right)  ^{n-1}\left(  -\sqrt{b^{2}-4c}\right)  \left(
-\frac{ax}{\sqrt{b^{2}-4c}}\right)  M_{n-1}^{\left(  2\right)  }\left(
\frac{ax}{\sqrt{b^{2}-4c}}-1;\frac{R^{+}}{R^{-}}\right)
\end{align*}
for $b>0$ and%
\begin{align*}
\Psi_{n}\left(  x;a,b,c\right)   &  =\left(  R^{+}\right)  ^{n}\left(
1-\frac{R^{-}}{R^{+}}\right)  \ \frac{ax}{\sqrt{b^{2}-4c}}M_{n-1}^{\left(
2\right)  }\left(  \frac{ax}{\sqrt{b^{2}-4c}}-1;\frac{R^{+}}{R^{-}}\right) \\
&  =\left(  R^{+}\right)  ^{n-1}\sqrt{b^{2}-4c}\frac{ax}{\sqrt{b^{2}-4c}%
}M_{n-1}^{\left(  2\right)  }\left(  \frac{ax}{\sqrt{b^{2}-4c}}-1;\frac{R^{+}%
}{R^{-}}\right)
\end{align*}
for $b<0.$
\end{proof}

\begin{remark}
If we define the inner product%
\begin{equation}
\left\langle f,g\right\rangle _{\beta}=%
{\displaystyle\sum\limits_{k=0}^{\infty}}
f(k)g(k)\gamma^{k}\frac{\Gamma\left(  \beta+k\right)  }{k!}, \label{IP4}%
\end{equation}
\medskip we have \cite{MR1149380}%
\[
\left\langle M_{n}^{\left(  \beta\right)  },M_{m}^{\left(  \beta\right)
}\right\rangle _{\beta}=\frac{n!\Gamma\left(  \beta+n\right)  }{\gamma
^{n}\left(  1-\gamma\right)  ^{\beta}}\delta_{n,m},\quad n,m\geq0.
\]
Therefore,%
\begin{gather*}
\left\langle M_{n}^{\left(  0\right)  },M_{m}^{\left(  0\right)
}\right\rangle _{0}=\left\langle \left(  1-\frac{1}{\gamma}\right)
\ xM_{n-1}^{\left(  2\right)  }(x-1;\gamma),\left(  1-\frac{1}{\gamma}\right)
\ xM_{m-1}^{(2)}(x-1;\gamma)\right\rangle _{0}\\
=\left(  1-\frac{1}{\gamma}\right)  ^{2}%
{\displaystyle\sum\limits_{k=1}^{\infty}}
k^{2}M_{n-1}^{\left(  2\right)  }(k-1;\gamma)M_{m-1}^{\left(  2\right)
}(k-1;\gamma)\gamma^{k}\frac{\Gamma\left(  k\right)  }{k!}\\
=\left(  1-\frac{1}{\gamma}\right)  ^{2}%
{\displaystyle\sum\limits_{k=0}^{\infty}}
M_{n-1}^{\left(  2\right)  }(k;\gamma)M_{m-1}^{\left(  2\right)  }%
(k;\gamma)\gamma^{k}\left(  k+1\right) \\
=\left(  1-\frac{1}{\gamma}\right)  ^{2}\left\langle M_{n-1}^{\left(
2\right)  },M_{m-1}^{\left(  2\right)  }\right\rangle _{2}\\
=\left(  1-\frac{1}{\gamma}\right)  ^{2}\frac{\left(  n-1\right)
!\Gamma\left(  n+1\right)  }{\gamma^{n-1}\left(  1-\gamma\right)  ^{2}}%
\delta_{n-1,m-1}=\frac{\left(  n-1\right)  !n!}{\gamma^{n+1}}\delta
_{n,m},\quad n,m\geq1.
\end{gather*}
Hence, $\left\{  M_{n}^{\left(  0\right)  }\ |\ n\geq1\right\}  $is an
orthogonal set with respect to the inner product $\left\langle \cdot
,\cdot\right\rangle _{0}$ defined by (\ref{IP4}).
\end{remark}

\begin{remark}
We will now show that $\left\{  M_{n}^{\left(  0\right)  }\ |\ n\geq0\right\}
$is an orthogonal set with respect to the inner product%
\[
\left\langle f,g\right\rangle =f(0)g(0)+%
{\displaystyle\sum\limits_{k=0}^{\infty}}
\Delta f(k)\Delta g(k)\frac{\gamma^{k}}{1-\gamma},
\]
where%
\[
\Delta f(x)=f\left(  x+1\right)  -f\left(  x\right)  .
\]

From (\ref{M0}) we have%
\begin{equation}
M_{n}^{\left(  0\right)  }(0;\gamma)=\delta_{n,0},\quad n\geq0, \label{M00}%
\end{equation}
which implies%
\begin{equation}
\left\langle M_{0}^{(0)},M_{n}^{(0)}\right\rangle =M_{n}^{(0)}(0;\gamma)+%
{\displaystyle\sum\limits_{k=0}^{\infty}}
0\frac{\gamma^{k}}{1-\gamma}=\delta_{n,0},\quad n\geq0. \label{IP0}%
\end{equation}
Using the formula \cite{MR1149380}
\[
\Delta M_{n}^{\left(  \beta\right)  }(x;\gamma)=n\left(  1-\frac{1}{\gamma
}\right)  M_{n-1}^{\left(  \beta+1\right)  }(x;\gamma),\quad n\geq0
\]
and (\ref{M00}), we get%
\begin{gather}
\left\langle M_{n}^{(0)},M_{m}^{(0)}\right\rangle =M_{n}^{(0)}(0;\gamma
)M_{m}^{(0)}(0;\gamma)\nonumber\\
+%
{\displaystyle\sum\limits_{k=0}^{\infty}}
n\left(  1-\frac{1}{\gamma}\right)  M_{n-1}^{\left(  1\right)  }%
(k;\gamma)m\left(  1-\frac{1}{\gamma}\right)  M_{m-1}^{\left(  1\right)
}(k;\gamma)\frac{\gamma^{k}}{1-\gamma}\label{IPM0}\\
=n^{2}\left(  1-\frac{1}{\gamma}\right)  ^{2}\frac{\left(  n-1\right)
!\Gamma\left(  n\right)  }{\gamma^{n-1}\left(  1-\gamma\right)  ^{2}}%
\delta_{n-1,m-1}\nonumber\\
=\frac{1}{\gamma^{n}}\left(  n!\right)  ^{2}\delta_{n,m},\quad n,m\geq
1.\nonumber
\end{gather}
Thus, (\ref{IP0}) and (\ref{IPM0}) give%
\[
\left\langle M_{n}^{(0)},M_{m}^{(0)}\right\rangle =\frac{1}{\gamma^{n}}\left(
n!\right)  ^{2}\delta_{n,m},\quad n,m\geq0.
\]

Some extensions using a similar inner-product were studied in \cite{MR1749794}%
, \cite{MR1752153}, \cite{MR1741786} and \cite{MR1884936}.
\end{remark}

\end{document}